\documentclass{conm-p-l}

\newtheorem{theorem}{Theorem}[section]
\newtheorem{lemma}[theorem]{Lemma}

\theoremstyle{definition}
\newtheorem{corollary}[theorem]{Corollary}

\newtheorem{proposition}[theorem]{Proposition}
\newtheorem{remark}[theorem]{Remark}

\theoremstyle{remark}

\newcommand{\Bp}{ {\mathcal D}^{+} }
\newcommand{\co}{ \mathcal O}
\newcommand{\cwe}{\mathbb C[w]^{(0)}}
\newcommand{\cwo}{\mathbb C[w]^{(1)}}
\newcommand{\C}{ {\mathbb C} }
\newcommand{\D}{ {\mathcal D} }
\newcommand{\DO}{ {\mathcal D}^{ { \mathcal O}} }
\newcommand{\dinf}{ d_{\infty} }
\newcommand{\Dinf}{ {\overline{d}}_{\infty} }

\newcommand{\Fl}{ {\mathcal F}^{\bigotimes l} }

\newcommand{\gl}{{gl} }
\newcommand{\glpm}{ \widehat{{ \mathfrak {gl} }}_{\pm} }

\newcommand{\hBp}{ \widehat{\mathcal D}^{+} }
\newcommand{\hD}{ \widehat{\mathcal D} }
\newcommand{\hf}{ \frac12}
\newcommand{\hgl}{\widehat{gl}_{\infty} }
\newcommand{\hL}{ \widehat{\Lambda} }
\newcommand{\hphi}{ \widehat{\phi}}

\newcommand{\vac}{ |0 \rangle }
\newcommand{\winf}{{\mathcal W}_{1 + \infty} }

\newcommand{\Z}{ {\mathbb Z} }

\numberwithin{equation}{section}

\begin{document}

\title{Dual pairs and infinite dimensional Lie algebras}

\author{Weiqiang Wang}
\address{Department of Mathematics, Yale University, 
New Haven, CT 06520}
\email{wqwang@math.yale.edu}


\begin{abstract}
  We construct and study various dual pairs between finite
dimensional classical Lie groups and infinite dimensional Lie algebras
in some Fock representations. 
The infinite dimensional Lie algebras here
can be either a completed infinite rank affine Lie algebra,
the $\winf$ algebra or its certain Lie subalgebras. 
We give a formulation in the framework of vertex algebras.  
We also formulate several conjectures and open problems.
\end{abstract}

\maketitle

\section*{Introduction}
 The theory of dual pairs of R.~Howe \cite{H1, H2} 
has many important applications to representation theory 
of reductive groups and $p$-adic groups etc. 
A dual pair is, roughly speaking, a pair of 
``maximally'' commuting  (reductive) Lie groups/algebras 
acting on a certain minimal module. The theory of dual pairs
also offers new insights
to various subjects such as invariant theory, spherical harmonics,
tensor product problems and so on. 

The aim of this paper is to
give a somewhat informal discussion of some of our recent
results \cite{W} \cite{KWY} (also see \cite{FKRW}) generalizing
dual pairs to an infinite dimensional setting 
and to formulate some conjectures and open problems.
The paper is mostly expository except the last section
which contains some new points.

A dual pair in our infinite dimensional setting
consists of a finite dimensional classical 
reductive Lie group and an infinite dimensional Lie algebra
acting on some Fock representations of free fermions or bosons.
In \cite{W}, we constructed
various dual pairs involving Lie groups $ O_n$, $Sp_{2n}$,
$Spin_n$ and Lie supergroup $Osp_{1, 2n}$ on one hand, and
Lie subalgebras of a completed infinite rank affine algebra 
$\hgl$ of $B, C, D$ types on the other hand.
We obtained isotypic decompositions with respect to dual pair 
actions and found explicit formulas for highest
weight vectors in all cases. In the hindsight, a duality of this sort
between a finite dimensional general linear group and $\hgl$
was first given in \cite{F2}, although heavy machinery rather
than the principles of dual pairs was used there. 

It turns out that dual pairs in our infinite dimensional setting
exhibit a new phenomenon
which is not seen in the finite dimensional cases
\cite{FKRW} \cite{KWY}: a natural Lie algebra may replace
another natural one acting on a same minimal module
with its dual pair partner fixed. Let us denote by $\hD$ the universal
central extension of the Lie algebra of differential operators
on the circle. In the literature $\hD$ is often denoted by
$\winf$ as well thanks to its connections to $\mathcal W$-algebras
\cite{BS, FeF}. There is an intimate relation between
modules of $\hgl$ and of $\hD$ \cite{KR} via a homomorphism
$\phi$ from $\hD$ to $\hgl$. It is realized \cite{FKRW}
that $GL_l$ and $\hD$ form a dual pair acting
on the Fock space of $l$ pairs of free fermions by replacing
$\hgl$ by $\hD$. While $\hgl$ is easier to
understand because of its similarity with 
finite dimensional general linear Lie algebras, $\hD$ is a more natural
object of study from the viewpoint of vertex algebras \cite{B, FLM, K}.

In view of the homomorphism $\phi$ from $\hD$
to $\hgl$, it is natural to study Lie subalgebras of
$\hD$ corresponding via $\phi$ to Lie subalgebras
of $\hgl$ of $B, C, D$ types. We found \cite{KWY}
there are two (families of) Lie subalgebras,
denoted by $\hD^{\pm}$, of $\hD$ fixed by two (families of)
anti-involutions of $\hD$ which do the job. 
In this way we obtain various dual pairs between classical Lie groups
and $\hD^{\pm}$ by substituting the Lie subalgebras 
of $\hgl$ of $B, C, D$ types appearing in the dual pairs
obtained in \cite{W} with appropriate $\hD^+$ or $\hD^-$.
The duality results formulated in terms of $\hD$ and
$\hD^{\pm}$ admit a natural reformulation in
the framework of vertex algebras. 

The length and complication of \cite{W} and \cite{KWY}
are to a large extent caused by the fact that we had to deal with
a dozen of different dual pairs.
In order to simplify matters and to focus on the main ideas, 
we restrict ourselves in this paper
to the consideration of dual pairs of $(A, A)$ and $(D, D)$ types in detail 
and comment on how other dual pairs can be constructed
while omitting almost all proofs. In the last section we work out some low
rank cases explicitly and show how these recapture and shed new light on 
some known results \cite{KRa} \cite{DG} \cite{DN}. 
We discuss several open problems and conjectures in the end.

This paper is an extended version of my talk in the conference
``Representations of affine and quantum affine algebras and
their applications'' held in NCSU at Raleigh
in May 1998. A similar talk was also presented in the Oberwolfach
meeting in July 1998. I would like to thank the organizers, Naihuan Jing
and Kailash Misra, respectively Richard Borcherds and Peter Slodowy, 
for the invitation.
\section{Lie algebras $\hgl$ and $ \dinf$}
\subsection{Lie algebra $\hgl$}
 We denote by $\gl$ the Lie algebra of all matrices
$(a_{ij})_{i,j \in \Z}$ such that $a_{ij} =0$
for $|i -j| $ large enough. 
Denote by $E_{ij}$ the infinite matrix with $1$ at the $(i, j)$-th
place and $0$ elsewhere. Let the weight of $ E_{ij}$ be $j - i $. This defines
a $\Z$--principal gradation $\gl = \bigoplus_{j \in \Z} \gl_j$.
Denote by $\hgl = \gl \bigoplus \C C$ the central extension
given by the following $2$--cocycle \cite{DJKM} \cite{KP}:
\begin{eqnarray}
 C(A, B) = \text{Tr} \,\left(
                     [J, A]B \right)
  \label{eq_cocy}
\end{eqnarray}
where $J = \sum_{ i \leq 0} E_{ii}$. 

The $\Z$--gradation
of Lie algebra $\gl$ extends to $\hgl$ by putting weight $C$ to be $ 0$. 
In particular, we have a triangular decomposition
$$\hgl = \widehat{\gl}_{+}  
          \bigoplus \widehat{\gl}_{0} \bigoplus \widehat{\gl}_{-}  $$
where 
$$ \glpm = \bigoplus_{ j \in \mathbb N}
\widehat{\gl}_{\pm j}, \quad 
\widehat{\gl}_{0} = \gl_0 \oplus \C C.$$
Denote by $\epsilon_i$ the linear
function on ${\mathfrak {gl} }_0$, s.t. 
$\epsilon_i (E_{jj} ) = \delta_{ij} (i, j \in \Z ).$
Then the root system of $\gl$ is 
$\Delta = \{ \epsilon_i - \epsilon_j, i, j \in \Z, i \neq j \}.$
The compact anti-involution $\omega$ is defined by
$\omega (E_{ij}) = E_{ji}$.

Given $c \in \mathbb C$ and $\Lambda $ in the restricted dual
${\gl}_{0}^{*}$ of ${\gl}_{0}$, we let
\begin{eqnarray*}
 {}^a\lambda_i & = & \Lambda (E_{ii}), \quad i \in \Z,    \\
 {}^a H_i     & = & E_{ii} - E_{i+1, i+1} + \delta_{i,0} C,  \\
  {}^a h_i & = & \Lambda( {}^a H_i ) = 
   \lambda_i - \lambda_{i+1} + \delta_{i,0} c.
 \label{eq_11}
\end{eqnarray*}
The superscript $a$ here denotes $\hgl$ which is of A type.

Denote by $L(\hgl; \Lambda, c)$ (or simply $L(\hgl; \Lambda)$ 
when the central charge is obvious from the text)
the $\hgl$--module with
highest weight $\Lambda$ and central
charge $c$. It is easy to see that $L(\hgl; \Lambda, c)$ is {\em quasifinite}
(namely having finite dimensional graded subspaces
according to the principal gradation of $\hgl$) if
and only if all but finitely many $h_i, i \in \Z$ are zero.
A quasifinite representation of $\hgl$ is {\em unitary}
if an Hermitian form naturally induced by $\omega$
is positive definite. 

Define ${}^a \Lambda_j \in \gl_0^* , j \in \Z$ as follows:
\begin{equation}
{}^a \Lambda_j ( E_{ii} ) =    
  \begin{cases}  
    1, & \text{for}\quad 0 < i \leq j \\
        -1, & \text{for}\quad j < i \leq 0  \\
        0, & \text{otherwise.}
  \end{cases}
\end{equation}
Define ${}^a \hL_0 \in \widehat{gl}_{0}^{*}$ by 
$$ {}^a \hL_0 (C) = 1, \quad
  {}^a \hL_0 (E_{ii}) = 0 \mbox{ for all } i \in \Z $$
and extend ${}^a \Lambda_j\/$ from $ \gl_0^{*}\/$ to 
$\widehat{gl}_{0}^{*}$ by letting
${}^a \Lambda_j (C) = 0.$ Then 
$${}^a \hL_j = {}^a \Lambda_j + {}^a \hL_0, \quad j \in \Z  $$ 
are the fundamental weights,
i.e. ${}^a \hL_j ( H_i ) = \delta_{ij}.$

One can prove that 
$L(\hgl; \Lambda, c)$ is unitary if and only if
$\Lambda =  {}^a \hL_{m_1} + \ldots + {}^a \hL_{m_k}$
for some $m_1 \ldots, m_k \in \Z$, and $ c = k \in \Z_{+}$.
\subsection{Lie algebra $\dinf$}
  \label{subsec_dinf}
   Let us consider the infinite dimensional 
vector space $\C [t, t^{-1}]$ and take a basis
$v_i = t^{ -i}, i \in \Z$. The Lie algebra $\gl$ acts on this vector
space naturally, namely $E_{ij} v_k = \delta_{jk}v_i$. We denote
by $\Dinf$ the Lie subalgebra of $\gl$ preserving the 
symmetric bilinear form 
$ D(v_i, v_j) = \delta_{i, 1-j},$ $i, j \in \Z.$
Namely we have 
$$\Dinf = \{ g \in \gl \mid D( a(u), v) + D(u, a(v)) = 0 \}
 = \{ (a_{ij})_{i,j \in \Z} \in \gl \mid a_{ij} = -a_{1-j,1-i} \}.$$

Denote by $\dinf = \Dinf \bigoplus \C C $ the central extension
given by the $2$-cocycle (\ref{eq_cocy}) restricted to $\Dinf$.
Then $\dinf$ has a natural triangular decomposition induced
from $\hgl$:
$\dinf = {\dinf}_{+} \bigoplus {\dinf}_0 \bigoplus {\dinf}_{-} .$

The set of simple coroots of $\dinf$, denoted by $ \Pi^\vee $
can be described as follows:
\begin{eqnarray*}
      \Pi^\vee & =& \{ \alpha_0^\vee = E_{00} + E_{-1,-1}
              -E_{2, 2} -E_{1, 1} + C,
                           \\
      & &\quad \alpha_i^\vee = E_{i, i} + E_{-i,-i}
        - E_{i+1,i+1} - E_{1-i,1-i}, i \in \mathbb N \}.
\end{eqnarray*}
             
Given $\Lambda \in {\dinf}_0^* $, we let
\begin{eqnarray*}
 {}^d \lambda_i & = &
        \Lambda (E_{ii} - E_{1-i, 1-i}) \quad (i \in \mathbb N), \\
 {}^d  H_i & = & E_{ii} + E_{-i, -i} - E_{i+1, i+1} - E_{-i+1, -i+1}
                                \quad (i \in \mathbb N),   \\
 {}^d  h & = & \Lambda( {}^d H_i ) = 
   \lambda_i  - \lambda_{i+1}  \quad (i \in \mathbb N), \\
 {}^d H_0 & = & E_{0,0} + E_{-1,-1} -E_{2,2} -E_{1,1} + 2C,   \\
 c & = & \hf ({}^d h_0 + {}^d h_1) + \sum_{i \geq 2} {}^d h_i.
 \label{eq_dcartan} 
\end{eqnarray*}
The superscript $d$ denotes $\dinf$ which is of D type.
Then we denote by ${}^d \hL_i $ the $i$-th fundamental weight
of $\dinf$, namely ${}^d \hL_i ({}^d H_j ) = \delta_{ij}.$
Denote by $L(\dinf; \Lambda, c)$ or simply $L(\dinf; \Lambda)$
the $\dinf$--module with highest weight $\Lambda$ and central charge $c$. 
We denote again by $\omega$ the anti-involution on $\dinf$
induced from $\omega $ on $\hgl$ by restriction.
A representation is {\em unitary}
if an Hermitian form naturally defined with respect to
the anti-involution $\omega$ on $\dinf$ is positive definite. 
$L(\dinf, \Lambda)$ is unitary if and only if
$\Lambda =  {}^d \hL_{m_1} + \ldots + {}^d \hL_{m_k}$, for some
$m_1, \ldots, m_k \in \Z_+$ and $c = \hf k \in \hf \Z_+.$
\section{Lie algebras $\hD$ and $\hD^{\pm}$}
 \subsection{Lie algebra $\D$}
 Let $\D_{as}$ be the associative algebra of differential operators on
the circle. $\D_{as}$ admits two distinguished bases: one is given
by 
\begin{eqnarray*}
  J^l_k = - t^{l+k} ( \partial_t )^l \quad
    ( l \in \Z_{+}, k \in \Z )
\end{eqnarray*} 
where $ \partial_t$ denotes $\frac{d}{dt}$; the other is given by
\begin{eqnarray*}
  L^l_k = - t^{k} D^l \quad
    (  l \in \Z_{+}, k \in \Z)   \nonumber
\end{eqnarray*}
where $D = t \partial_t$. It is easy to see that $J^l_k = -t^k [D]_l$.
Here and further we use the notation 
$
 [x]_l = x(x-1)\ldots (x-l+1). 
   \label{eq_symb}
$
                     
Let $\D$ denote the Lie algebra obtained from $\D_{as}$ by
taking the usual Lie bracket.
Denote by $\hD$ the central extension of $ \D $
by a one-dimensional center with a generator $C$: $\hD = \D + \C C$.
The Lie algebra $\hD$
has the following commutation relations \cite{KR}
\begin{eqnarray*}
  \left[
     t^r f(D), t^s g(D)
  \right]
    & = & t^{r+s} 
    \left(
      f(D + s) g(D) - f(D) g(D+r) 
    \right) \nonumber \\ 
   & + & \Psi 
      \left(
        t^r f(D), t^s g(D)
      \right)
      C
\end{eqnarray*}
where 
\begin{eqnarray}
  \Psi 
      \left(
        t^r f(D), t^s g(D)
      \right)   = 
   \begin{cases}
        \sum_{-r \leq j \leq -1} f(j) g(j+r),& r= -s \geq 0  \\
        0, & r + s \neq 0. 
   \end{cases}
  \label{eq_cocykp}
\end{eqnarray}

Set the weight of $J^l_k$ to be $ k$ and the weight of $ C$  $0$.
This defines the principal $\Z$-gradations of $\D_{as}, \D$ and $\hD$:
\begin{eqnarray*}
  \D = \bigoplus_{j \in \Z} {\mathcal D}_j, \quad
  \hD = \bigoplus_{j \in \Z} \widehat{\mathcal D}_j .
\end{eqnarray*}
Thus we have the triangular decomposition
\begin{eqnarray*}
  \hD = \hD_{+} \bigoplus \hD_{0}  \bigoplus \hD_{-},
\end{eqnarray*}
where 
\begin{eqnarray*}
  \hD_{\pm} = \bigoplus_{j \in \pm \mathbb N} \widehat{\mathcal D}_j,
  \quad
  \hD_{0} = {\D}_{0}  \bigoplus {\mathbb C} C.
\end{eqnarray*}

By abuse of notation, we use again $J^l_k, L^l_k$ to denote
the corresponding elements in $\hD$.
  The elements $J^0_m, m\in \Z$ span a Heisenberg algebra:
\[
 [ J^0_m , J^0_n ] = \delta_{m, -n} m C.
\]
The elements $J^1_m, m\in \Z$ span a Virasoro algebra:
\[
 [ J^1_m , J^1_n ] = (m -n) J^1_{m +n} - \delta_{m , -n} \frac{m^3 - m}6  C.
\]

\subsection{Lie algebras $\hD^{\pm}$}
 Recall that an {\em anti-involution} 
$\sigma$ of $\D_{as}$ is an involutive anti-automorphism of $\D$, i.e.
$ \sigma^2 = I, \sigma (a X + bY) =a\sigma (X) +b \sigma (Y)$ and
$ \sigma (XY)=\sigma (Y)\sigma (X),$ 
where $a, b \in \mathbb C, X, Y \in \D.$

\begin{proposition}
  Any anti-involution $\sigma$ of $\D_{as}$ preserving the
 principal $\Z$-gradation is one of the following:
 \begin{enumerate}
  \item[1)] $\sigma_{-,b}( t) = -t, \quad \sigma_{-,b}(D) = -D +b;$ 
  \item[2)] $\sigma_{+,b}( t)= t, \quad \sigma_{+,b}(D) = -D +b, 
     \quad b \in \C. $
 \end{enumerate}
\end{proposition}

It follows immediately that
 $
  \sigma_{\pm, b}  (\partial_t) = \mp \left(\partial_t - (b+1) t^{-1} \right).
 $ 
Given $ s \in \C$, denote by $\Theta_s$ the automorphism of $\D_{as}$ which sends
$t$ to $t$ and $D$ to $D + s$. Equivalently $\Theta_s$ is
given by sending $a \in \D$ to $t^{-s}a t^s $, 
the conjugate of $a$ by $t^s$.
Clearly $\Theta_s$
preserves the principal $\Z$-gradation of $\D_{as}$.
We have 
\[
  \sigma_{\pm, b} \circ \Theta_s = \sigma_{\pm, b + s}, \quad
  \Theta_{-s} \circ \sigma_{\pm, b} = \sigma_{\pm, b + s} . 
\]
                  
Denote by ${\mathcal D}^{\pm, b}$ the Lie subalgebra 
of $\D$ fixed by $- \sigma_{\pm, b}$, namely 
$$ {\mathcal D}^{\pm, b} = \{ a \in \D \mid \sigma_{\pm, b} (a) = -a \}. $$
Since $\sigma_{\pm, b}$ preserves the principal $\Z$-gradation of $\D$
there is an induced $\Z$-gradation on $ {\mathcal D}^{\pm, b}$:
$ {\mathcal D}^{\pm, b} = \bigoplus_{j \in \mathbb Z} {\mathcal D}^{\pm, b}_j,$
where $ {\mathcal D}^{\pm, b}_j = {\mathcal D}^{\pm, b}\cap \hD_j .$
The relation among $ {\mathcal D}^{\pm, b}$ for different $b \in \C$ 
can be described as follows.
\begin{lemma}
  The Lie algebras $ {\mathcal D}^{+, b}$ (resp. $ {\mathcal D}^{-, b}$) 
 for different $b \in \C$ 
  are all isomorphic. More precisely, we have
  $ \Theta_s ( {\mathcal D}^{\pm, b}) = {\mathcal D}^{\pm, b -2s}.$
    \label{lem_iso}
\end{lemma}
 
Thanks to Lemma \ref{lem_iso} we may fix $b$ in $ {\mathcal D}^{\pm, b}$.
It turns out there is a best choice by putting $b = -1$
(partly because $\sigma_{\pm, -1} ( \partial_t) = \mp \partial_t$).
From now on we set $ \Bp = {\mathcal D}^{+, -1}$ and $\D^- = {\mathcal D}^{-, -1}$
(note that our $\D^-$ is different from the one in \cite{KWY}).
Accordingly we set $\D^{\pm}_j = \D^{\pm, -1}_j$ and so
$\D^{\pm} = \sum_{j \in \Z} \D^{\pm}_j$.
We shall see that there is a canonical vertex algebra structure
on the vacuum module of the central extension of $\Bp$.

Let us denote by $\cwe$ (resp. $\cwo$) the set of all even 
(resp. odd) polynomials in $ \C [w]$. We let $ \overline{k} = 0$
if $k$ is an odd integer and $ \overline{k} = 1$ if $k$ is even.
\begin{lemma}
    We have
 \begin{eqnarray*}
   {\mathcal D}^{+}_j & = & \left\{ t^j g \left(D + (j +1)/2 \right)
      \mid  g(w) \in \cwo \right \},  \\
   {\mathcal D}^{-}_j & = & \left\{ t^j g \left( D + (j +1)/2 \right)
      \mid  g(w) \in \C [w]^{(\bar{j})} \right \} , \; j \in \Z.
 \end{eqnarray*}
  \label{lem_graded}
\end{lemma}

We will concentrate on the Lie algebra $\D^+$. The story for
$\D^-$ is quite similar in any event. We set
\[
W^{n}_k = - \hf t^k \left(
                    [D ]_n - [ -D -k -1 ]_n
                \right) \quad ( k \in \Z, n \in \mathbb N ).
\]
One can check that $ W^{n}_k$ $( k \in \Z, n \in {\mathbb N}_{odd}
=\{1, 3, 5, \ldots \})$ form a basis of $\D^+$.
                  
By abuse of notation we again denote by $\Psi$ the restriction of the
2-cocycle $\Psi$ to $\Bp$.
Denote by $\hBp$ the central extension of $ \Bp $
by $ \C C$ corresponding to $\Psi$. The Lie algebra
$\hBp$ is a subalgebra of $\hD$ by definition.
Note that $W^1_k \; (k \in \Z)$ span a Virasoro algebra, namely we have
$$ \left[ W^1_m, W^1_n 
  \right]
  = (m -n) W^1_{m +n} + \delta_{m,-n} \frac{m^3 -m}{12} C.
$$

\begin{remark}  \rm  \label{rem_twist}
  We observe from Lemma~\ref{lem_graded} that $\hD^+_j$
 coincides with $\hD^-_j$ for $j$ even.
 The relation between $\hD^+$ and $\hD^-$ is analogous to
 the relation between untwisted affine algebras and
 twisted affine algebras (it is more than an analog!).
\end{remark}
\section{Quasifinite modules of $\hD$ and $\hD^+$ and a link to $\hgl$}
\subsection{Quasifinite modules of $\hD$ and $\hD^+$}
   Let $ \mathfrak g =\oplus_{j \in \mathbb Z} \mathfrak g_j$ 
be a $\Z$-graded Lie algebra satisfying
$ [{\mathfrak g}_i, {\mathfrak g}_j] \subset {\mathfrak g}_{i+j}.$ Let
${\mathfrak g}_{\pm} =\oplus_{j > 0} {\mathfrak g}_{ \pm j}$.
A $\mathfrak g$-module $V$ is called {\it $\Z$-graded} if 
$ V=\oplus_{j \in \Z} V_j$ and ${\mathfrak g}_i V_j \subset V_{j -i}. $
A graded $\mathfrak g$-module is called {\it quasifinite} if 
$ \dim V_j < \infty$ for all $j$. 

Since $\mathfrak g$ carries a triangular decomposition
${\mathfrak g} = {\mathfrak g}_+ \oplus {\mathfrak g}_0 \oplus {\mathfrak g}_-$,
we can construct as usual a Verma module $ M(\mathfrak g; \Lambda)$
with a highest weight vector annihilated by ${\mathfrak g}_+$
and highest weight $\Lambda \in {\mathfrak g}_0^*$.
Denote by $L(\mathfrak g; \Lambda)$ 
the irreducible quotient of $M(\mathfrak g; \Lambda)$ 
by the maximal proper graded submodule.

We specialize $\mathfrak g$ to $\hD$ and $\hBp$ from now on. 
Since the Cartan subalgebras of $\hD$ and $\hBp$ are infinite
dimensional, the classifications of quasifinite highest
weight $\hD$-modules and $\hBp$-modules turns out
to be quite non-trivial.

For a given $\Lambda \in (\D_0)^*$, 
we characterize $\Lambda$ by its {\em labels} $\Delta_n = -\Lambda( D^n)$ 
$(n \in \Z_+)$ and the central charge $c=\Lambda (C)$. 
Similarly for a given $\Lambda \in (\D_0^+)^*$, 
we characterize $\Lambda$ by its labels $\Delta^+_n = -\Lambda( (D +1/2)^n)$ 
$(n \in \mathbb N_{{odd}} )$
and the central charge $c=\Lambda (C)$. 

We introduce the generating functions
\[
 \Delta_{\Lambda}(x) = \sum_{n \geq 0 } 
 \frac{x^n}{n!} \Delta_n, \quad
 \Delta^+_{\Lambda}(x) = \sum_{n \in {\mathbb N}_{{ odd} } } 
 \frac{x^n}{n!} \Delta^+_n.
\]
We call $\Delta_{\Lambda}(x)$, resp. $\Delta^+_{\Lambda}(x)$,
the {\em characteristic generating function} of $L(\hD; \Lambda)$,
resp. of $L(\hD^+; \Lambda)$.
Sometimes we simply drop the subscript ${\Lambda}$. 
Clearly we have
\begin{eqnarray}
 \Delta (x) & = & - \Lambda \left( e^{xD} \right), \\
 \Delta^+ (x) & = & - \hf \Lambda \left( e^{x(D + 1/2)} - e^{ -x(D + 1/2)} \right).
   \label{eq_triv}
\end{eqnarray}

A {\it quasipolynomial} is a finite
linear combination of functions of the form 
$p(x) e^{ \alpha x}$, where $p(x)$ is a polynomial 
and $\alpha \in \C$. The following beautiful characterization of 
quasi-finiteness of an irreducible module $ L(\hD; \Lambda)$
is due to Kac and Radul \cite{KR}.
\begin{theorem}
 A $\hD$-module $L(\hD; \Lambda)$ is quasifinite if and only if
 $$
   \Delta (x) = \frac{ F(x)}{e^x -1}
 $$
 where $F(x)$ is a quasipolynomial such that $F(0) = 0$.
   \label{th_kr}
\end{theorem}

Given an irreducible quasifinite highest weight $\hD$-module $V$ with
central charge $c \in \C$ and with characteristic generating function
$\Delta (x)$, we write $F(x) + c$ as a finite sum of the form
$
  \sum_i p_i (x) e^{s_i x},
$
where $p_i (x)$ are non-zero polynomials satisfying
$\sum_i p_i (0) = c $
and $s_i$ are distinct complex numbers.

We call $s_i$ the {\em exponents} of $V$ with {\em multiplicities}
$p_i (x)$. We single out a class of important 
quasifinite $\hD$-modules called
{\em positive primitive} modules.
$V$ is called {\em positive primitive}
if all the multiplicities are positive integers. We observe that
$V$ is uniquely determined by its exponent-multiplicity set ${\mathcal E}$.
We will denote $V$ by $L( \hD; {\mathcal E})$.

The following characterization of 
quasi-finiteness of an irreducible highest weight $\hD^+$-module
is obtained in \cite{KWY}.
\begin{theorem}
 A $\hD^+$-module $L(\hD^+; \Lambda)$ is quasifinite if and only if
 $$
   \Delta^+ (x) = \frac{ F(x)}{2 \sinh (x/2 )}
 $$
 where $F(x)$ is an even quasipolynomial such that $F(0) = 0$.
   \label{th_delta}
\end{theorem}

Given an irreducible quasifinite highest weight $\hBp$-module $V$ with
central charge $c \in \C$ and with characteristic generating function
$\Delta^+ (x)$, we write $F(x) + c$ as a finite sum of the form
 \begin{eqnarray}
  \sum_i p_i (x) \cosh (e_i^+ x) + \sum_j q_j (x) \sinh (e_j^- x),
   \label{eq_ka}
 \end{eqnarray}
where $p_i (x)$ (resp. $q_j (x)$) are non-zero even (resp. odd) polynomials
and $e_i^+$ (resp. $ e_j^-$) are distinct complex
numbers. Clearly $\sum_i p_i (0) = c. $ The expression (\ref{eq_ka})
is unique up to a sign of $e_i^+$ or a simultaneous change of signs of
$ e_j^-$ and $q_j (x)$.

We call $e_i^+$ (resp. $e_j^-$) the {\em even type} (resp. {\em odd type})
{\em exponents} of $V$ with {\em multiplicities}
$ p_i (x)$ (resp. $ q_j (x)$). 
$V$ is called {\em positive primitive}, if there is
no odd type exponents and the multiplicities
of its (even type) exponents $e_i^+$ are integers $n_i$ such that
$n_i >0$ when $e_i \neq  0$ and
 $- \hf n_{i_0} \leq n_i \in \hf \Z$ when $e_i =  0$. Here
 ${i_0}$ is the index such that $e_{i_0} = 1$.
We denote by $\mathcal E$ the set of (even type) nonzero exponents
 with their multiplicities for a positive primitive module $V$.
 The pair $({\mathcal E}, c)$ determines uniquely the module $V$. We will
 denote this primitive module by $L( \hBp; {\mathcal E}, c)$.
\subsection{A homomorphism from $\hD$ to $\hgl$}
  
Let $\co $ be the algebra of all holomorphic functions on $\C$ 
with the topology of uniform convergence on compact sets. Define
$
 {\mathcal O}^o  =  \{ f \in {\co} \mid f(w) = - f(-w) \}.
$
We define a completion $\D^{\mathcal O}$ of $\D$ consisting  of 
all differential operators of the form $ t^j f(D)$ where
$f \in \co$ and $j \in \Z$. This induces a completion 
$\D^{+, \mathcal O}$ of $\Bp$.
Many things defined for $\D$ and $\D^+$, such as $2$-cocycle
$\Psi$, commutation relations etc, extend naturally to their
completions.

Recall that the Lie algebra $\gl$ acts on the vector
space $\C [t, t^{-1}]$ by $E_{ij} v_k = \delta_{jk} v_i$,
where $v_k = t^{-k}$. Note that the Lie algebras
$\D$ and $\D^{\mathcal O}$ naturally acts on the same space
as differential operators:
$
 t^j f(D) t^k = f(k)t^{k +j}.
$
In this way we obtain an algebra embedding
$\phi$ of $\D$ into $\gl$ given by
$
 \phi \left( t^k f( D) \right) = \sum_{j \in \Z} f (-j )  E_{j-k, j}.
$
We remark that $\phi$ extending to $\DO$ becomes a surjective homomorphism. 

When restricted to $\D^+$ and  $\D^{+, \mathcal O}$, we have
\begin{eqnarray*}
  \phi \left( t^k f( D + (k +1)/2) 
          \right) = \sum_{j\in \Z}
      f(-j +(k +1)/2 )  E_{j-k, j}
\end{eqnarray*}
for $ f \in {\mathcal O}^o.$ The image of 
$\phi: \D^{+, \mathcal O} \rightarrow \gl$
can be easily shown to be $\Dinf$. 

It turns out that the $2$-cocycle $\Psi$ (\ref{eq_cocykp}) on $\D$ and $\D^+$
and the $2$-cocycle given by (\ref{eq_cocy}) on $\gl$ are compatible via
the homomorphism $\phi :\D \rightarrow \gl$. Thus we can
extend $\phi$ to 
a homomorphism $ \hphi$ from $\hD$ to $\hgl$ by letting
\begin{equation}
 \hphi \left( t^k f( D) \right) = \sum_{j \in \Z} f (-j )  E_{j-k, j},
  \quad  \hphi (C)  =  C.
  \label{eq_embd}
\end{equation}
When restricted to $\hD^+$, we obtain a homomorphism
$\hphi : \hD^+ \rightarrow \dinf.$

Note that the $\Z$-gradations of $\hD$ and $\hD^+$ are compatible
with those on $\hgl$ and $\dinf$ respectively via $\hphi$. 
Thus we may regard $\hD$ (resp. $\hD^+$) as a {\em dense} subalgebra of
$\hgl$ (resp. $\dinf$) via $\hphi$.
\begin{proposition}
 The pullback of an irreducible quasifinite module over $\hgl$ (resp. $\dinf$)
 via the homomorphism $\hphi$ is an irreducible quasifinite
 module over $\hD$ and (resp. $\hD^+$).
\end{proposition}
\section{(A, A) and (D, D) duality}
  \subsection{Some actions on a Fock space $\Fl$}
  Let us take $l$ pairs of free fermionic fields
$$ 
 \psi^{\pm p}(z) = \sum_{n \in \hf +\Z}
  \psi^{\pm p}_n z^{ -n -\hf }.
$$
The anti-commutation
relations among $\psi^{\pm p}_n$ are equivalent to the following
operator product expansions 
$$ \psi^{+ p}(z) \psi^{- q}(w) \sim \frac{\delta_{p q}}{z-w}, \quad
   \psi^{+ p}(z) \psi^{+ q}(w) \sim 0, \quad
   \psi^{- p}(z) \psi^{- q}(w) \sim 0,
$$
where $p,q = 1, \ldots, l$. We denote by $\Fl$
the Fock space generated by a vacuum vector $\vac$ annihilated by
$\psi^{\pm p}_n$, $n > 0$, $p = 1, \ldots, l$.

Introduce the following generating functions
\begin{eqnarray}
  E(z,w) & \equiv & \sum_{i,j \in \Z} E_{ij} z^{i-1 }
     w^{-j  } 
   = \sum_{p =1}^l :\psi^{+,p} (z) \psi^{-,p} (w):  \label{eq_genelin} \\
  e^{pq} (z) & \equiv & \sum_{n \in \Z} e^{pq} (n) z^{ -n-1 }
    =  : \psi^{-,p}(z) \psi^{-,q} (z):   \nonumber    \\ 
  e_{**}^{pq} (z) & \equiv & \sum_{n \in \Z}
                            e_{**}^{pq}(n)z^{ -n-1 }
    =  : \psi^{+,p}(z) \psi^{+,q} (z): \nonumber     \\
  e_*^{pq} (z) & \equiv & \sum_{n \in \Z} e_*^{pq} (n) z^{ -n-1 }
    =  : \psi^{+,p}(z) \psi^{-,q} (z): 
         \nonumber 
\end{eqnarray}
where $p,q = 1, \dots, l$, and the normal ordering $::$ means that 
the operators annihilating $\vac$ are moved to the right and 
multiplied by $-1$.

It is well known \cite{F1, KP} that the operators $e^{pq}(n), e_*^{pq}(n), 
e_{**}^{pq}(n)$, $p,q = 1, \dots, l$, $n \in \Z$ defined as above
form a representation of the affine algebra
$\widehat{so}_{2l}$ of level $ 1$. In particular
$e_*^{pq}(n)$, $p,q = 1, \dots, l$, $n \in \Z$ 
form a representation of the affine algebra
$\widehat{gl}_l$ of level $ 1$.
The operators $e^{pq}(0), e_*^{pq}(0), e_{**}^{pq}(0)
 \;(p, q = 1, \cdots, l)$ form
the horizontal subalgebra $\mathfrak {so}_{2l}$ in $\widehat{ \mathfrak {so}}_{2l}$
while $e_*^{pq}(0) \;(p, q = 1, \cdots, l)$ form
the horizontal subalgebra $\mathfrak {gl}_l$ in $\widehat{ \mathfrak {gl}}_l$.
We identify the Borel subalgebra ${\mathfrak b}( {\mathfrak {so}}_{2l} )$
with the one generated by $e_{**}^{pq} \;(p \neq q), e_*^{pq}\; (p \leq q),
p, q = 1, \cdots, l.$ We take the Borel subalgebra of $gl_l$ to be
${\mathfrak b}( {\mathfrak {gl}}_l )= \mathfrak {gl}_l \cap
{\mathfrak b}( {\mathfrak {so}}_{2l} )$.

It follows from (\ref{eq_genelin}) that
 \begin{eqnarray}
    &&  \sum_{i,j \in \Z} ( E_{ij} - E_{1-j,1-i} ) z^{i- 1} w^{-j} \\
   &=  &
   \sum^l_{k=1} (: \psi^{+,k} (z) \psi^{-,k} (w) : 
    - : \psi^{+,k} (w) \psi^{-,k} (z) :).  \nonumber
   \label{eq_gener}
 \end{eqnarray}
One can show \cite{DJKM}
that $ E_{ij}( i,j \in \Z)$ defined above span $\hgl$ and therefore
$ E_{ij} - E_{1-j, 1-i} \;( i,j \in \Z)$
span $\dinf$ with central charge $l$. 
\subsection{Duality involving $\hgl$ and its subalgebras}
   We observe that the horizontal subalgebra ${\mathfrak {gl}}_l$ 
acting on $\Fl$ lifts to $GL_l$. It is not difficult to show that
the action of $\hgl$ commutes with the action of $GL_l$ 
on $\Fl$. We can show that $\hgl$ and $GL_l$ form a dual pair acting on $\Fl$
by a similar argument as in the finite dimensional case \cite{H1}.

It is well known that the set $\Sigma (A)$ of
the $n$-tuples $(m_1, \cdots, m_l),$
$ m_1 \geq m_2 \geq \ldots \geq m_l $, $m_i \in \Z$
parametrizes the finite dimensional irreducible $GL_l$-modules.
We define a map $\Lambda^{ \mathfrak{aa}}$ 
from $ \Sigma (A)$ to $\widehat{gl}_0^*$ as follows:
given $ \lambda = (m_1, \cdots, m_l),$ we define 
\[
 \Lambda^{ \mathfrak{aa}} (\lambda ) = 
  \sum_{i =1}^l {}^a \hL_{m_i}.
\]
The following duality of $(A, A)$ type
is due to I.~Frenkel \cite{F2} (also see \cite{FKRW}).

\begin{theorem}  [Frenkel]
  Under the joint action of $GL_l$ and $\hgl$, we have the
   isotypic decomposition
   \begin{eqnarray*}
     \Fl = \bigoplus_{\lambda \in \Sigma (A) } 
          V(GL_l; \lambda) \otimes L \left(\hgl;
                                 \Lambda^{\mathfrak{aa}} (\lambda), l
                               \right)
   \end{eqnarray*}
   where 
   $ L \left(\hgl; \Lambda^{\mathfrak{aa}} (\lambda), l  \right)$
   is the irreducible highest weight $\hgl$-module of highest weight
   $ \Lambda^{\mathfrak{aa}} (\lambda)$ and central charge $l$.
  \label{th_IF}
\end{theorem}

In \cite{F2} heavy machinery rather than principles of dual pairs
was used in establishing the theorem.
An alternative simple proof goes as follows, see \cite{W} for detail.
It follows directly from the principles of dual pairs
that the decomposition is multiplicity-free. We find 
explicit highest weight vectors and read off the highest
weights with respect to the action of $GL_l$ and $\hgl$ 
respectively. The theorem is now established by the simple observation
that all highest weights for finite dimensional
irreducible $GL_l$-modules already show up.

\begin{remark}  \rm
  All unitary irreducible highest weight $\hgl$-modules with central charge $l$
 appear in the above decomposition. They are in one-to-one
 correspondence with all irreducible finite dimensional $GL_l$-modules.
\end{remark}

We note that the horizontal subalgebra ${\mathfrak {so}}_{2l}$
acting on $\Fl$ lifts to $SO_{2l}$ and extends to an
action of $O_{2l}$ naturally. One can show (see \cite{W}) that
the actions of $O_{2l}$ and $\dinf$ on $\Fl$ commute
and moreover $\dinf$ and $ O_{2l}$ form a dual pair.

$O(2l)$ is a semi-direct product of $SO(2l)$ by $\Z /2 \Z$.
Below we will use a highest weight to denote the corresponding
module for short. If $\lambda$ is a representation of $SO(2l)$ of
highest weight $(m_1, m_2, \ldots, m_l )$ ($m_l \neq 0$),
then the induced representation of $\lambda$ 
to $O(2l)$ is irreducible and its restriction to
$SO(2l)$ is a sum of $(m_1, m_2, \ldots, m_l )$ and
$(m_1, m_2, \ldots, - m_l )$.
We denote this irreducible representation of $O(2l)$ by
$(m_1, m_2, \ldots, \overline{m}_l )$, where $m_l$ is chosen to be 
greater than 0. If $ m_l = 0$, the representation
$\lambda = (m_1, m_2, \ldots, m_{l-1}, 0 )$ extends to two
different representations of $O(2l)$, i.e. 
$\lambda$ itself and $\lambda \bigotimes det$,
where $det$ is the $1$-dimensional non-trivial representation of $O(2l)$. 

Thus the irreducible modules of $O_{2l}$ are parametrized by the set
\begin{eqnarray*}
  \Sigma(D) =
   & & \left\{ (m_1, m_2, \ldots, \overline{m}_l ) \mid
       m_1 \geq m_2 \geq \ldots \geq m_l > 0, m_i \in \Z;  \right.
                               \\
   & &   (m_1, m_2, \ldots, m_{l-1}, 0 ) \bigotimes det, \\
   & &   \left.   (m_1, m_2, \ldots, m_{l-1}, 0 ) \mid
       m_1 \geq m_2 \geq \ldots \geq  m_{l-1} \geq 0, m_i \in \Z \right\}. 
\end{eqnarray*}
We denote by $V(O_{2l}; \lambda)$ the irreducible
$O_{2l}$-module corrsponding to $\lambda \in \Sigma(D)$ from now on.

  We define a map 
$\Lambda^{ \mathfrak{dd}}: \Sigma (D) \longrightarrow {\dinf}_0^* $
by sending $ \lambda = (m_1, \cdots, \overline{ m}_l) \; (m_l > 0 )$ to 
$$  \Lambda^{ \mathfrak{dd}} (\lambda) =
   (l -i ) \; {}^d \hL_0 +
         ( l -i )\; {}^d \hL_1 + \sum_{k =1}^i {}^d \hL_{m_k}, 
$$
sending $(m_1, \cdots, m_j, 0, \ldots, 0 )\; (j <l )$ to 
$$  \Lambda^{ \mathfrak{dd}} (\lambda) =
   (2l -i -j)\; {}^d \hL_0
        + (j-i ) \;{}^d \hL_1 + \sum_{k =1}^i {}^d \hL_{m_k}, 
$$
and sending $ (m_1, \ldots, m_j, 0, \ldots, 0 ) \bigotimes {det} 
   \; ( j < l)$ to 
 $$  \Lambda^{\mathfrak{dd}} (\lambda) =
  (j-i)\; {}^d \hL_0 
      + (2l -i -j)\; {}^d \hL_1 + \sum_{k =1}^i {}^d \hL_{m_k}, $$
if $m_1 \geq \ldots m_{i} > m_{i+1} = \ldots = m_{j} =1
   > m_{j+1} = \ldots = m_l = 0.$

We have the following duality of $(D, D)$ type.

\begin{theorem} [\cite{W} ]
  Under the joint action of $O_{2l}$ and $\dinf $, we have the
   isotypic decomposition
 \begin{eqnarray*}
   \Fl = \bigoplus_{\lambda \in \Sigma (D) } 
          V(O_{2l}; \lambda) \otimes L \left(\dinf;
                                 \Lambda^{\mathfrak{dd}} (\lambda), l
                               \right)
 \end{eqnarray*}
 where $L \left(\dinf; \Lambda^{\mathfrak{dd}} (\lambda), l  \right)$
 is the irreducible $\dinf$-module of highest weight
 $ \Lambda^{\mathfrak{dd}} (\lambda)$ and central charge $l$.
  \label{th_Fpair}
\end{theorem}
\begin{remark}  \rm
  All unitary irreducible highest weight $\dinf$-modules 
 with central charge $l$
 appear in the above decomposition. They are in one-to-one
 correspondence with all irreducible finite dimensional $O_{2l}$-modules.
\end{remark}

\begin{remark} \rm   \label{rem_vari}
 We indicate below that there are many variations of Theorems~\ref{th_IF} and 
 Theorem~\ref{th_Fpair} which give rise to different dual pairs
 involving with various Lie groups. See \cite{W} for
 more detail.
 \begin{enumerate}
  \item If we consider the Fock space of $l$ pairs of fermions 
 together with
 a neutral fermion in a similar way, we can construct a dual pair
 between $O_{2l +1}$ and $\dinf$ (with central charge $l + \hf$).

  \item If we change the indices for the Fourier components
 of the fermions $\psi^{\pm p} (z)$ $(p =1, \ldots , l)$
 from half-odd-integers to integers,
 we can construct a dual pair between the pin group $Pin_{2l}$ and 
 a Lie subalgebra of $\hgl$ of type $B$
 acting on the Fock space modified accordingly. If in addition we add
 a neutral fermion as in the previous case, the spin group
 $Spin (2l +1)$ will show up instead.

 \item If we consider the Fock space of free bosons instead of free fermions,
 we obtain dual pairs involving the symplectic
 group $Sp_{2n}$ and Lie supergroup $Osp_{1, 2n}$.

 \item It turns out that the finite dimensional Lie groups arising in the dual
 pairs considered so far always arise from lifting of a horizontal
 subalgebra of an appropriate untwisted affine algebra acting on
 a Fock space \cite{FF}. If we consider the horizontal subalgebras
 of certain twisted affine algebras acting on Fock spaces, we obtain 
 other new dual pairs between finite dimensional Lie groups and
 infinite dimensional Lie agebras as well.
 \end{enumerate}
\end{remark}
\section{Lie algebras $\hD$, $\hD^+$ and vertex algebras}
\subsection{Duality involving $\hD$ and $\hD^+$}
  We obtain an action of $\hD$ on $\Fl$ by composing
the action of $\hgl$ and the homomorphism $\hphi$
given by the formula (\ref{eq_embd}). It follows that
the actions of $GL_l$ and $\hD$ on $\Fl$ commute and that
$GL_l$ and $\hD$ form a dual pair on $\Fl$.

Introduce the generating function
\[
  J^n (z) = \sum_{k \in \Z} J^n_k z^{ -k -n -1}.
\]
\begin{lemma}  \label{lem_simple}
   On $\Fl$ we have
 \[
 J^n (z) = \sum^l_{k=1} : \partial_z^n \psi^{-,k} (z) \psi^{+,k} (z) : .
 \]
\end{lemma}

For $ \lambda = (m_1, \cdots, m_l)\in \Sigma(A) $ where
$m_1 \geq \ldots \geq m_l$, $m_i \in \Z$,
we let ${\mathcal E} (\lambda)$ be the set of exponents 
$m_k$ $(k =1, \ldots, l)$ with multiplicity $1$. Note that
the $\hD$-module $L(\hD; {\mathcal E} (\lambda))$ is positive primitive.

\begin{remark}  \rm
  We implicitly assumed that $m_1, \ldots, m_l$ are distinct
in defining ${\mathcal E} (\lambda)$ above. In general
if $m', \ldots, m_{\alpha}'$ are distinct numbers among
$m_1, \ldots, m_l$ with multiplicities $a_1, \ldots, a_{\alpha}$,
then ${\mathcal E} (\lambda)$ should be the set of exponents 
$m', \ldots, m_{\alpha}'$ with multiplicities $a_1, \ldots, a_{\alpha}$.
With this understood,
we have prefered the oversimplified and slightly wrong way of
defining ${\mathcal E} (\lambda)$ here and below.
\end{remark}

\begin{theorem} [\cite{FKRW}]  \label{th_fkrw}
 Under the joint action of $GL_l$ and $\hD$,
 we have the isotypic decomposition
 \[
   \Fl = \bigoplus_{\lambda \in \Sigma (A)} 
          V( GL_l; \lambda) \otimes L \left(\hD;
                                 {\mathcal E} (\lambda)
                               \right).
 \]
\end{theorem}

\begin{corollary}   \label{cor_glinv}
   The space of invariants of $GL_l$ in the Fock space $\Fl$ is 
 naturally isomorphic to the irreducible 
 module $ L (\hgl; l\; {}^a \hL_0 )$,
 or equivalently to the irreducible module 
 $L \left(\hD; {\mathcal E}(0) \right)$.
\end{corollary}

We obtain the action of $\hBp$ on $\Fl$ by composing
the action of $\dinf$ and the homomorphism $\hphi$
given by the formula (\ref{eq_embd}). Similarly 
$O_{2l}$ and $\dinf$ form a dual pair. Introduce
the generating function
\[
  W^n (z) = \sum_{k \in \Z} W^n_k z^{ -k -n -1}.
\]
\begin{lemma}   \label{lem_ffr}
   On $\Fl$ we have
 \[
 W^n (z) = \hf \sum^l_{k=1} (: \partial_z^n \psi^{-,k} (z) \psi^{+,k} (z) : 
    + : \partial_z^n \psi^{+,k} (z) \psi^{-,k} (z) :).
 \]
\end{lemma}

For $ \lambda = (m_1, \cdots, \overline{ m}_l)\in \Sigma(D) $ where
$m_1 \geq \ldots \geq m_i > m_{i +1} = \ldots = m_l = 1$,
we let ${\mathcal E} (\lambda)$ be the set of exponents $m_k$ $(k =1, \ldots, i)$
with multiplicity $1$ and the exponent $1$ with multiplicity $l -i$;
for $\lambda = (m_1, \cdots, m_l)\in \Sigma(D)$ where 
$$
  m_1 \geq \ldots m_{i} > m_{i+1} = \ldots = m_{j} =1
   > m_{j+1} = \ldots = m_l = 0 \;( j < l),
$$ 
we let ${\mathcal E} (\lambda)$ be the set of 
exponents $m_k$ $(k =1, \ldots, i)$ of multiplicity $1$,
exponent $1$ of multiplicity $j -i$;
for $ (m_1, \cdots, m_l) \bigotimes {det} \in \Sigma(D)$ where 
$$
 m_1 \geq \ldots m_{i} > m_{i+1} = \ldots = m_{j} =1
   > m_{j+1} = \ldots = m_l = 0\; ( j < l),
$$
we let ${\mathcal E} (\lambda)$ be the set of 
exponents $m_k$ $(k =1, \ldots, i)$ of multiplicity $1$,
exponent $1$ of multiplicity $2l -i -j$.
We will write $(0, \cdots, 0 ) $
and $(0, \cdots, 0 ) \bigotimes {det}$ for short as $0$ and $det$ respectively.
Note that the $\hD^+$-module $L(\hD^+ ; {\mathcal E} (\lambda), l)$
is positive primitive.
\begin{theorem} [\cite{KWY}]   \label{th_kwy}
  Under the joint action of $O_{2l}$ and $\hD^+$,
 we have the isotypic decomposition
 \begin{equation}
   \Fl = \bigoplus_{\lambda \in \Sigma (D)} 
          V( O_{2l}; \lambda) \otimes L \left(\hBp;
                                 {\mathcal E} (\lambda), l
                               \right).
 \end{equation}
  \label{th_Wpairdd}
\end{theorem}

\begin{corollary}
   The space of invariants of $O_{2l}$ in the Fock space $\Fl$ is 
 naturally isomorphic to the irreducible 
 module $ L (\dinf; 2l\; {}^d \hL_0 )$ of central charge $l$,
 or equivalently to the irreducible module 
 $L \left(\hBp; {\mathcal E}(0), l \right)$.
   \label{cor_dl}
\end{corollary}
          
  The Dynkin diagram of $\dinf$ admits an automorphism of order 2 denoted
 by $\sigma$. $\sigma$ induces naturally an automorphism of order 2
 of $\dinf$, which is denoted again by $\sigma$.
 $\sigma$ acts on the set of highest weights for $\dinf$
 by mapping $ \lambda = \; {}^d h_0 \; {}^d  \hL_0 + \; {}^d h_1 \;{}^d \hL_1 
  + \sum_{i \geq 2} \; {}^d h_i  \; {}^d \hL_i$ to
 $\sigma(\lambda ) = \; {}^d h_1  \; {}^d \hL_0 + \; {}^d h_0  \; {}^d \hL_1 
  + \sum_{i \geq 2} \; {}^d h_i  \; {}^d \hL_i$. 
 In this way one obtain an irreducible module of the semi-product
 $\sigma \propto \dinf$
 on $ L(\dinf; \lambda) \bigoplus L(\dinf; \sigma(\lambda ) )$
 if $\sigma(\lambda ) \neq \lambda $ and on $ L(\dinf; \lambda) $
 if $\sigma(\lambda ) = \lambda $.  Note that 
 $( SO_{2l}, \sigma \propto \dinf)$ form a dual pair on $\Fl$. 

\begin{corollary}
 The space of invariants of $SO_{2l}$ in $\Fl$ is 
 isomorphic to the $\dinf$-module
 $L(\dinf; 2l \; {}^d \hL_0 ) \bigoplus L(\dinf; 2l \; {}^d \hL_1 )$
 or equivalently the $\hBp$-module $L \left(\hBp; {\mathcal E}(0), l \right)
 \bigoplus L \left(\hBp; {\mathcal E} (det), l \right)$. 
\end{corollary}

\begin{remark}  \rm
 The dual pairs as indicated in Remark~\ref{rem_vari}
 can be reformulated in terms of $\hD^+$ or $\hD^-$ accordingly
 as in Theorem~\ref{th_kwy}.
\end{remark}
\subsection{Vertex algebras associated to $\hD$ and $\hD^+$}
  The irreducible quasifinite $\hD$-module
 $L (\hD; {\mathcal E} (0) )$ 
 (resp. $\hD^+$-module $L (\hD^+; {\mathcal E} (0), c)$)
is usually refered to as the irreducible
vacuum $\hD$-module (resp. $\hD^+$-module), denoted by $V_c$
(resp. $V_c^+$). Denote by $\vac$ the highest weight vector.

The fermionic Fock space $\Fl$ is one of the simplest
examples of vertex algebras \cite{B} \cite{FLM} \cite{K}.
One may view $GL_l$ and $O_{2l}$ as automorphism groups
of the vertex algebra $\Fl$. The space
$V_l$, being isomorphic to the space of invariants of
$GL_l$ of the vertex algebra $\Fl$ by Corollary~\ref{cor_glinv}, 
inherits a vertex algebra 
structure. Similarly it follows from Corollary~\ref{cor_dl}
that $V_l^+$ is a vertex subalgebra of $\Fl$.

We want to show that $V_c$ and $V_c^+$ are vertex algebras
for general $c$. One can easily check that
\begin{eqnarray*}
\lim_{z \rightarrow 0} J^l (z) \vac = J^l_{-l -1} \vac, &&
 \quad  [ J^1_0, J^l (z) ] = \partial J^l (z), \\
\lim_{z \rightarrow 0} W^l (z) \vac = W^l_{-l -1} \vac, &&
 \quad  [ W^1_0, W^l (z) ] = \partial W^l (z).
\end{eqnarray*}

For a typical element in $V_c$ (resp. $V_c^+$) of the form
$  J^{n_j}_{-k_j -n_j -1} \cdots J^{n_1}_{-k_1 -n_1 -1} \vac$
(resp. $ W^{n_j}_{-k_j -n_j -1} \cdots W^{n_1}_{-k_1 -n_1 -1} \vac$)
where $ k_i \in \Z_+$,
we associated a field (state-field correspondence)
$:\partial_z^{(k_j)} W^{n_j}(z) \cdots \partial_z^{(k_1)} W^{n_1}(z): .$

One can prove by using Lemmas~\ref{lem_simple} and 
\ref{lem_ffr} that
\[
  (z -w)^{m +n +2} \left[ J^m (z), J^n (w) \right] = 0, \quad
  (z -w)^{m +n +2} \left[ W^m (z), W^n (w) \right] = 0.
\]
in $\Fl$ for all $l \in \Z_+$. It follows 
that these equations hold for any $c$.

Combining all the above, we have proved the following 
theorem \cite{FKRW} \cite{KWY}. The simplicity here follows 
from the irreducibility of $V_c$ (resp. $V_c^+$) as a module over
$\hD$ (resp. $\hD^+$).
\begin{theorem}
  $V_c$ and $V_c^+$ are simple vertex algebras.
\end{theorem}

\begin{corollary}  \label{cor_module}
 $L(\hD; {\mathcal E} (\lambda) ),$ $ \lambda \in \Sigma (A)$
 are irreducible modules of the vertex algebra $V_l$;
 $L(\hD^+; {\mathcal E} (\lambda), l ),
 $ $ \lambda \in \Sigma (D)$ are irreducible
 modules of $V^+_l$.
\end{corollary}

\begin{remark} \rm
 We obtain a dual pair between a Lie group and
 a vertex algebra $(GL_l, V_l)$  (resp. $(O_{2l}, V^+_l)$)
 by substituting $ \hD$ (resp. $ O_{2l}$) with $V_l$ (resp. $V^+_l$).
\end{remark}
\section{Some examples and open problems}
\subsection{Examples}
  Let us first consider the case $l =1$. 
Set ${\mathcal F} = {\mathcal F}^{\bigotimes 1}$ and 
$\psi^{\pm} (z) = \psi^{\pm, 1} (z)$. We claim that
the $(GL_1, \hD)$-duality in Theorem~\ref{th_fkrw}
is essentially the celebrated boson-fermion correspondence.
To see that, we notice an infinitesimal generator of $GL_1$,
$J^0_0 =\sum_{n \in \Z} :\psi^-_{-n}\psi^+_n :$, acts on a typical element
$\psi^-_{-m_1} \ldots \psi^-_{-m_j}
\psi^+_{-n_1} \ldots \psi^+_{-n_i} \vac$ as multiplication by $j -i$.

Recall that the Fourier components of
$J^0(z)$ generate a Heisenberg algebra $\mathcal H$. One 
easily shows that $J^n (z) ( n \geq 1)$ can be written as
a normally ordered polynomial in terms of $ \partial^i J^0 (z)$,
$i \in \Z^+$. For example, 
$J^1 (z) = \hf (\partial  J^0 (z) + : J^0 (z) J^0 (z):)$.
Thus an irreducible representation of $\hD$ inside $\mathcal F$ 
remains irreducible when restricted to the Heisenberg algebra
$\mathcal H$. Theorem~\ref{th_fkrw} for $l =1$ says that
\[
 {\mathcal F} = \bigoplus_{n \in \Z} {\mathcal H}^n
\]
where ${\mathcal H}^n$, the subspace of $\mathcal F$ on
which $J^0_0$ acts as multiplication by $n$, is an irreducible module
over $\mathcal H$. In particular the vertex algebra
$V_1$ is isomorphic to the Heisenberg vertex algebra ${\mathcal H}^0$
generated by $J^0 (z)$. 

Now consider the $(O_2, \hD^+)$-duality on ${\mathcal F}$
in Theorem~\ref{th_kwy}. $O_2$ is the semi-direct product of
$SO_2$ by $\Z / 2\Z$ generated by the order $2$ automorphism
$\tau $ which exchange
$\psi^+$ and $\psi^-$. An irreducible $O_2$-module
except the trivial module and the determinant module (which
are of $1$-dimension) is two-dimensional. When restricted 
to  $SO_2$, it decomposes into 
a sum of two $SO_2$-modules, with weight say $m$ and $-m$, 
$m \in \mathbb N$. 

Theorem~\ref{th_kwy} (or rather together with Theorem~{3.2}, \cite{W}) 
has the following consequences:
${\mathcal H}^m $ $( m \neq 0)$ is irreducible as a module over $ \hD^+$
and ${\mathcal H}^m $ is isomorphic to ${\mathcal H}^{ -m} $
as $ \hD^+$-modules; ${\mathcal H}^0$ is a direct sum of 
the $(\pm 1)$-eigenspaces ${\mathcal H}_{\pm}^0$ of $\tau$;
${\mathcal H}_{+}^0$ is isomorphic to 
$L(\hD^+; {\mathcal E}(0), 1)$ (which is in turn isomorphic to $V_1^+ $)
while ${\mathcal H}_-^0$ is isomorphic to $L(\hD^+; {\mathcal E}(det), 1)$.

The $q$-character formulas for irreducible
$\hD^+$-modules (counting the graded dimensions) appearing in 
$\Fl$ were obtained in \cite{KWY}. In particular
one deduces that
\[
 ch_q V_1^+
      =
        \frac1{ 2 \prod_{j \geq 1} (1 -q^{2j})}
       \left( \prod_{j \geq 1} (1 +q^{j}) +\prod_{j \geq 1} (1 -q^{j})
       \right) 
      = \frac{\sum_{j \geq 0} (-q)^{j^2} }{ \prod_{j \geq 1} (1 -q^{j})}.
\]
The second equlity here is obtained by using the Gauss identity
\[
 \frac{\prod_{j \geq 1} (1 -q^{j})}{\prod_{j \geq 1} (1 +q^{j}) }
 = \sum_{j \in \Z} (-q)^{j^2} .
\]

Denote by $L(c, h)$ the irreducible module of Virasoro algebra 
with highest weight $h$ and central charge $c$. 
It is well known (cf. e.g. \cite{KRa}) that the character of $L(1, 4m^2)$ is 
\[
 ch_q L(1, 4m^2) =
 \frac{q^{4m^2} - q^{(2m +1)^2} }{ \prod_{j \geq 1} (1 -q^{j})} .
\]
It follows from the $q$-character formula of $V_1^+$ that 
$V_1^+ = \bigoplus_{ m \geq 0} L(1, 4m^2)$ 
as a module over the Virasoro algebra generated by $W^1 (z)$. 
Similarly we can show that
\[
  L(\hD^+; {\mathcal E}(det), 1) = \bigoplus_{ m \geq 0} L(1, (2m +1)^2)
\]
 as a module over this Virasoro algebra.

One can show that $W^n (z)$ for all $n$ can be written
as normally ordered polynomials in terms of
$ \partial^i W^1 (z)$ and $ \partial^i W^3 (z)$, $i \in \Z_+$.
One checks that $W^1 (z)$ and $W^3 (z)$
generate the $W(2, 4)$ algebra with central 
charge $1$ \cite{BS}. This shows that
$V_1^+$ is isomorphic to the $W(2, 4)$ algebra. The latter is
simple as we know that $V_c^+$ is simple for any $c$.
 
\begin{remark}  \rm
The fixed-point vertex subalgebra ${\mathcal H}_+^0$ 
of ${\mathcal H}^0$ by $\tau$ was studied in detail 
in \cite{DG} \cite{DN},
where ${\mathcal H}_+^0$ is denoted by $M(1)^+$.
The results derived above mainly as consequences of 
the $(O_2, \hD^+)$-duality have been obtained
in \cite{DG, DN} in a very different way. 
Irreducible modules over the vertex algebra ${\mathcal H}_+^0$
has also been classified \cite{DN}. We note that this classification
can be achieved by the approach of \cite{W2} as well. In view of
the isomorphism between $V_1^+ $ and ${\mathcal H}_+^0$, the classification
of the irreducible modules over $V_1^+$ is done. Relations (conjecturally
being minimal)  among the two generators 
in $V_1^+$ could be obtained by our approach.
\end{remark}

Now let us consider the Fock space, denoted by 
${\mathcal F}^{\bigotimes \hf}$, of a neutral fermion 
$\varphi (z) = \sum_{ n \in \hf + \Z} \varphi_n z^{ -n -\hf}$
satisfying $ [\varphi_m, \varphi_n]_+ = \delta_{m ,-n}$.
One has a duality between $O_1 $ which is $\Z / 2 \Z$ and
$\dinf$ (or $\hD^+$) with central charge $\hf$ (see Remark~\ref{rem_vari} 
or \cite{W, KWY} for more detail). Thus the decomposition of 
${\mathcal F}^{\bigotimes \hf}$ into even and odd parts gives rise to
\[ 
 {\mathcal F}^{\bigotimes \hf} = L(\dinf; {}^d \hL_0 ) \bigoplus
  L(\dinf; {}^d \hL_1 ) ,
\]
or equivalently,
 \[ 
 {\mathcal F}^{\bigotimes \hf} = L(\hD^+; {\mathcal E}(0), 1/2 ) \bigoplus
  L(\hD^+; {\mathcal E}(det), 1/2 ) .
\]

In this case, we have $W^i(z) = \partial^i \varphi (z) \varphi (z) $.
The Fourier components of $W^1 (z)$ generate the
Virasoro algebra with central charge $\hf$. 
One can easily show by induction that
$W^n(z) $ for all $n$ can be
written as a normally ordered polynomial in terms of $W^1(z)$. 
It follows that $L(\hD^+; {\mathcal E}(0), \hf )$ (which 
is the even part of ${\mathcal F}^{\bigotimes \hf}$) is irreducible 
over the Virasoro algebra generated by $W^1 (z)$ ($c = \hf$)
with highest weight $0$, and $  L(\hD^+; {\mathcal E}(det), \hf )$ (which 
is the odd part of ${\mathcal F}^{\bigotimes \hf}$) is irreducible 
over this Virasoro algebra with highest weight  $\hf$.
This identification of the even (resp. odd) part of
${\mathcal F}^{\bigotimes \hf}$ with the irreducible module over
the Virasoro algebra was proved (cf. e.g. \cite{KRa}) by using
deep results on the structures of modules over the Virasoro algebra.
In particular we have proved that the vertex algebra
$V_{\hf}^+$ is isomorphic to the Virasoro vertex algebra
with central charge $\hf$. 
This is indeed the only $c$ such that $V_c^+$ is rational.

\begin{remark}  \rm
 We can obtain an explicit decomposition of a tensor product 
 of two minimal modules
 of Virasoro algebra with central charge $\hf$ into irreducible
 modules of Virasoro algebra with central charge $1$ from the
 duality results discussed above by using the see-saw pair technique
 in the theory of dual pairs.
\end{remark}
\subsection{Some open problems}
    We discuss several open problems and conjectures to conclude this paper.

{\bf Problem 1}: determine a minimal set of generators
and relations for the vertex algebras $V_c$ and $V_c^+$;
classify all the irreducible modules of $V_c$ and $V_c^+$.

Recall that one can associate a $\mathcal W$-algebra 
${\mathcal W} {\mathfrak g}$ 
to an arbitrary complex simple Lie algebra $\mathfrak g$ \cite{BS} \cite{FeF}. 
In general such a $\mathcal W$-algebra can
be defined in terms of intersections of screening operators,
cf. \cite{FeF}. With some mild modification, one can define
a $\mathcal W$-algebra ${\mathcal W}({\mathfrak gl}_l)$ associated to
${\mathfrak gl}_l$ \cite{FKRW}. 

 When $c \notin \Z$, the structure of $V_c$ is well understood.
In \cite{FKRW}, it is shown that the vertex algebra
$V_l$ associated to $\hD$ for $l \in \mathbb N$ is isomorphic
to ${\mathcal W}({\mathfrak gl}_l)$ with central charge $l$.
More precisely, the first $l$ generating fields $J^i (z)$,
$i =0, 1, \ldots, l -1$ freely generate $V_l$ while
the other  fields $J^i (z)$, $i  \geq l$ can be written as
some normally ordered polynomials in terms of  
$\partial^n J^i (z)$, $n \geq 0$,
$i =0, 1, \ldots, l -1$. The irreducible
modules of $V_l$ were classified.

\begin{remark}  \rm
The vertex algebra $V_c$ can be viewed as some universal
$\mathcal W$-algebra which is sort of 
``limit'' of the ${\mathcal W}({\mathfrak gl}_l)$ algebra
which are generated by fields of conformal dimension 
$1, 2, \ldots, l$. One naturally asks a similar
question for other series of $\mathcal W$ algebras
associated to classical simple Lie algebras such as
 ${\mathcal W} {\mathcal D}_l$. 
Since the  ${\mathcal W} {\mathcal D}_l$
algebra has generating fields of conformal dimension
$2, 4, \ldots, 2l-2$ and $l$ (which equals
$1$ plus the exponents of $D_l$), one cannot take the naive ``limit''
due to the existence
of the field with conformal dimension $l$.
However $V_c^+$ can be viewed as a universal $\mathcal W$-algebra 
which is sort of ``limit'' of ${\mathcal W} {\mathcal D}_l$
in the following sense \cite{KWY}: a direct sum
of $V_l$ with the distinguished irreducible module
$L(\hD^+; {\mathcal E}(det), l )$ over $V_l$ is isomorphic
to ${\mathcal W} {\mathcal D}_l$.
\end{remark}

The structure of $V_c$ for a negative integral $c$ 
is much more complicated. In \cite{W1}, we 
showed that $V_{-1}$ is isomorphic to the ${\mathcal W}({\mathfrak gl}_3)$
algebra. The irreducible modules of $V_{-1}$ were classified in \cite{W2, W1}.
A relation among the generating fields of $V_{-1}$ was implicit there
and conjecturally it determines all other relations. 
For a general negative integral $c = -l$ it is conjectured that
 a minimal set of generating fields in $V_{-l}$ consists of
the first $l^2 +2l$ fields $J^i (z)$, $ 0 \leq i < l^2 +2l$. 
The conjecture was verified in the 
case $c = -1$ \cite{W2}.

For $c$ generic which means $c \notin \hf \Z$ in this case, the
structure of $V_c^+$ is understood \cite{KWY}.
For $c \in \hf \Z_+$ 
we conjecture that a minimal set of generating fields of $V_c^+$
consists of the first $2c$ fields $W^i (z)$, 
$ i =1, 3, \ldots, 4c -1$. The conjecture in the case $c = \hf, 1$
was verified in the preceding subsection.

{\bf Problem 2}: clarify precise
relations between $\hD^+$ and $\hD^-$. 

As remarked in
Remark~\ref{rem_twist}, $\hD^-$ should be thought as
some kind of twisted algebra of $\hD^+$. 
We conjecture that irreducible
modules appearing in the duality decomposition involving
$\hD^-$ with central charge $c \in \hf \Z_+$ (see \cite{KWY})
can be viewed as irreducible twisted modules over the vertex algebra
$V_c^+$. A natural problem then will be to classify all the
irreducible twisted modules over $V_c^+$.

{\bf Problem 3}: formulate results in \cite{KWY}
in terms of lattice vertex algebras by means of the boson-fermion
correspondence.

\bibliographystyle{amsalpha}

\end{document}